\newif\ifconf                           
\newif\ifshort                          
\newif\ifcompress                       

\conffalse\shortfalse\compressfalse

\documentclass[12pt,twoside]{article}
\usepackage{latexsym,amsmath,amssymb,hyperref}
\def\eqvsp{}  \newdimen\paravsp  \paravsp=1.3ex
\def\,{\mskip 3mu} \def\>{\mskip 4mu plus 2mu minus 4mu} \def\;{\mskip 5mu plus 5mu} \def\!{\mskip-3mu}
\def\dispmuskip{\thinmuskip= 3mu plus 0mu minus 2mu \medmuskip=  4mu plus 2mu minus 2mu \thickmuskip=5mu plus 5mu minus 2mu}
\def\textmuskip{\thinmuskip= 0mu                    \medmuskip=  1mu plus 1mu minus 1mu \thickmuskip=2mu plus 3mu minus 1mu}
\textmuskip
\def\beq{\eqvsp\dispmuskip\begin{equation}}    \def\eeq{\eqvsp\end{equation}\textmuskip}
\def\beqn{\eqvsp\dispmuskip\begin{displaymath}}\def\eeqn{\eqvsp\end{displaymath}\textmuskip}
\def\bqa{\eqvsp\dispmuskip\begin{eqnarray}}    \def\eqa{\eqvsp\end{eqnarray}\textmuskip}
\def\bqan{\eqvsp\dispmuskip\begin{eqnarray*}}  \def\eqan{\eqvsp\end{eqnarray*}\textmuskip}

\newtheorem{theorem}{Theorem}
\newtheorem{corollary}[theorem]{Corollary}
\newtheorem{lemma}[theorem]{Lemma}
\newtheorem{definition}[theorem]{Definition}
\newenvironment{keywords}{\centerline{\bf\small
Keywords}\begin{quote}\small}{\par\end{quote}\vskip 1ex}

\def\paradot#1{\vspace{\paravsp plus 0.5\paravsp minus 0.5\paravsp}\noindent{\bf\boldmath{#1.}}}

\def\req#1{(\ref{#1})}
\def\toinfty#1{\smash{\stackrel{#1\to\infty}{\longrightarrow}}}
\def\eps{\varepsilon}
\def\nq{\hspace{-1em}}
\def\qed{\hspace*{\fill}\rule{1.4ex}{1.4ex}$\quad$\\}
\def\fr#1#2{{\textstyle{#1\over#2}}}
\def\SetR{I\!\!R}
\def\qmbox#1{{\quad\mbox{#1}\quad}}
\def\E{{\bf E}}

\def\lb{{\log_2}}
\def\g{\gamma}
\def\t{\theta}

\def\l{\ell}
\def\o{\omega}
\def\O{\Omega}
\def\Oo{{\Omega^\circ}}
\def\Ov{{\vec\Omega}}
\def\B{{\cal B}}
\def\X{{\cal X}}
\def\Y{{\cal Y}}
\def\F{{\cal F}}
\def\M{{\cal M}}
\def\P{{P}} 
\def\Prob{{\text{Pr}}}
\def\MDL{{\text{MDL}}}
\def\MDLI{{\text{MDLI}}}
\def\Bayes{{\text{Bayes}}}

\begin{document}

\title{\vspace{-4ex}
\vskip 2mm\bf\Large\hrule height5pt \vskip 4mm
Discrete MDL Predicts in Total Variation
\vskip 4mm \hrule height2pt}
\author{{\bf Marcus Hutter}\\[3mm]
\normalsize RSISE$\,$@$\,$ANU and SML$\,$@$\,$NICTA \\
\normalsize Canberra, ACT, 0200, Australia \\
\normalsize \texttt{marcus@hutter1.net \ \  www.hutter1.net}
}
\date{September 2009}
\maketitle

\begin{abstract}
The Minimum Description Length (MDL) principle selects the model
that has the shortest code for data plus model. We show that for a
countable class of models, MDL predictions are close to the true
distribution in a strong sense. The result is completely general. No
independence, ergodicity, stationarity, identifiability, or other
assumption on the model class need to be made. More formally, we
show that for any countable class of models, the distributions
selected by MDL (or MAP) asymptotically predict (merge with) the
true measure in the class in total variation distance. Implications
for non-i.i.d.\ domains like time-series forecasting, discriminative
learning, and reinforcement learning are discussed.
\def\contentsname{\centering\normalsize Contents}
{\parskip=-2.7ex\tableofcontents}
\end{abstract}

\vspace*{-2ex}

\section{Introduction}\label{secIntro}

The {\em minimum description length} (MDL) principle recommends to
use, among competing models, the one that allows to compress the
data+model most \cite{Gruenwald:07book}. The better the compression, the
more regularity has been detected, hence the better will predictions
be. The MDL principle can be regarded as a formalization of Ockham's
razor, which says to select the simplest model consistent with the
data.

\paradot{Multistep lookahead sequential prediction}
We consider sequential prediction problems, i.e.\ having {\em observed
sequence} $x\equiv (x_1,x_2,...,x_\l) \equiv x_{1:\l}$, {\em predict}
$z\equiv (x_{\l+1},...,x_{\l+h}) \equiv x_{\l+1:\l+h}$, then observe
$x_{\l+1}\in\X$ for $\l\equiv\l(x)=0,1,2,...$. Classical prediction is
concerned with $h=1$, multi-step lookahead with $1<h<\infty$, and
total prediction with $h=\infty$. In this paper we consider the
last, hardest case. An infamous problem in this category is the
Black raven paradox \cite{Maher:04,Hutter:07uspx}: Having
observed $\l$ black ravens, what is the likelihood that {\em all}
ravens are black. A more computer science problem is (infinite horizon)
reinforcement learning, where predicting the infinite future is
necessary for evaluating a policy. See Section \ref{secAppl} for
these and other applications.

\paradot{Discrete MDL and Bayes}
Let $\M=\{Q_1,Q_2,...\}$ be a {\em countable class of
models}=\linebreak[1]theories=\linebreak[1]hypotheses=\linebreak[1]probabilities
over sequences $\X^\infty$, sorted w.r.t.\ to their {\em
complexity=codelength} $K(Q_i)=2\lb i$ (say), containing the {\em
unknown true sampling distribution} $P$. Our main result will be for
arbitrary measurable spaces $\X$, but to keep things simple in the
introduction, let us illustrate MDL for finite $\X$.

In this case, we define $Q_i(x)$ as the $Q_i$-probability of data
sequence $x\in\X^\l$. It is possible to code $x$ in $\log
P(x)^{-1}$ bits, e.g.\ by using Huffman coding. Since $x$ is
sampled from $P$, this code is optimal (shortest among all prefix
codes).
Since we do not know $P$, we could select the $Q\in\M$
that leads to the shortest code on the observed data $x$. In order
to be able to reconstruct $x$ from the code we need to know which
$Q$ has been chosen, so we also need to code $Q$, which takes $K(Q)$
bits. Hence $x$ can be coded in $\min_{Q\in\M}\{-\log Q(x)+K(Q)\}$
bits. MDL selects as model the minimizer
\beqn
  \MDL^x \;:=\; \arg\min_{Q\in\M}\{-\log Q(x)+K(Q)\}
\eeqn

Given $x$, the true predictive probability
of $z$ is $P(z|x)=P(xz)/P(x)$. Since $P$ is unknown we use
$\MDL^x(z|x):=\MDL^x(xz)/\MDL^x(x)$ as a substitute.
Our main concern is how close is the latter to the former.
We can measure the distance between two predictive distributions
by
\beq\label{dh}
  d_h(P,Q|x) \;=\; \sum_{z\in\X^h}\big|P(z|x)-Q(z|x)\big|
\eeq
for $h<\infty$ and $d_\infty=\lim_{h\to\infty} d_h=\sup\{d_1,d_2,...\}$.
It is easy to see that $d_h$ is monotone increasing
and that $d_\infty$ is twice the total variation distance (tvd)
defined in \req{tvd}.

MDL is closely related to Bayesian prediction, so a comparison to
existing results for Bayes is interesting. Bayesians use $\Bayes(z|x)$ for
prediction, where $\Bayes(x):=\sum_{Q\in\M} Q(x)w_Q$ is the Bayesian
mixture with prior weights $w_Q>0\,\forall Q\in\M$ and
$\sum_{Q\in\M}w_Q=1$. A natural choice is $w_Q\propto 2^{-K(Q)}$.

\paradot{Results}
The following results can be shown
\beq\label{results}
  {\sum_{\l=0}^\infty \E[d_h(P,\MDL^x|x_{1:\l})] \leq 21\,h\!\cdot\!2^{K(P)},\quad \atop
  \sum_{\l=0}^\infty \E[d_h(P,\Bayes|x_{1:\l})] \leq h\!\cdot\!\ln w_P^{-1},\qquad }
  {d_\infty(P,\MDL^x|x) \to 0 \atop
  d_\infty(P,\Bayes|x)\to 0 }  \;\;\left\{
  {\mbox{almost surely} \atop
   \mbox{for}\;\; \l(x)\!\to\!\infty} \right.
\eeq
where the expectation $\E$ is w.r.t.\ $P[\cdot|x]$.
The left statements for $h<\infty$ imply $d_h\to 0$ almost surely,
including some form of convergence rate. For Bayes it has been proven
in \cite{Hutter:03spupper}; for MDL
the proof in \cite{Hutter:05mdl2px}
can be adapted.
As far as asymptotics is concerned, the right results $d_\infty\to
0$ are much stronger, and require more sophisticated proof
techniques. For Bayes, the result follows from \cite{Blackwell:62}.
The proof for MDL is the primary novel contribution of this
paper; more precisely for arbitrary measurable $\X$ in total
variation distance.
Another general consistency result is presented in
\cite[Thm.5.1]{Gruenwald:07book}. Consistency is shown (only) in
probability and the predictive implications of the result are
unclear.
A stronger almost sure result is alluded to, but the given reference
to \cite{Barron:91} contains only results for i.i.d.\ sequences
which do not generalize to arbitrary classes.
So existing results for discrete MDL are far less satisfactory than
the elegant Bayesian prediction in tvd.

\paradot{Motivation}
The results above hold for completely arbitrary countable model
classes $\M$. No independence, ergodicity, stationarity, identifiability, or other assumption need to be made.

The bulk of previous results for MDL are for continuous model
classes \cite{Gruenwald:07book}. Much has been shown for classes of
independent identically distributed (i.i.d.) random variables
\cite{Barron:91,Gruenwald:07book}. Many results naturally generalize
to stationary-ergodic sequences like ($k$th-order) Markov. For
instance, asymptotic consistency has been shown in \cite{Barron:85}.
There are many applications violating these assumptions, some of
them are presented below and in Section \ref{secAppl}.

One can often hear the exaggerated claim that (e.g.\ unlike Bayes)
MDL can be used even if the true distribution $P$ is not in $\M$.
Indeed, it can be used, but the question is wether this is any good.
There are some results supporting this claim, e.g.\ if $P$ is in the
closure of $\M$, but similar results exist for Bayes. Essentially
$P$ needs to be at least close to some $Q\in\M$ for MDL to work, and
there are interesting environments that are not even close to being
stationary-ergodic or i.i.d.

Non-i.i.d.\ data is pervasive \cite{Amini:09}; it includes all
time-series prediction problems like weather forecasting and stock
market prediction \cite{Cesa:06}. Indeed, these are also perfect
examples of non-ergodic processes. Too much green house gases, a
massive volcanic eruption, an asteroid impact, or another world war
could change the climate/economy irreversibly.
Life is also not ergodic; one inattentive second in a car can have
irreversible consequences. Also stationarity is easily violated in
multi-agent scenarios: An environment which itself contains a
learning agent is non-stationary (during the relevant learning
phase). Extensive games and multi-agent reinforcement learning are
classical examples \cite{Weinberg:04}.

Often it is assumed that the true distribution can be uniquely
identified asymptotically. For non-ergodic environments, asymptotic
distinguishability can depend on the realized observations, which
prevent a prior reduction or partitioning of $\M$. Even if
principally possible, it can be practically burdensome to do so,
e.g.\ in the presence of approximate symmetries. Indeed this problem
is the primary reason for considering {\em predictive} MDL.
MDL might never identify the true distribution, but our main result
shows that the sequentially selected models become predictively
indistinguishable.

The countability of $\M$ is the severest restriction of our result.
Nevertheless the countable case is useful. A semi-parametric problem
class $\bigcup_{d=1}^\infty\M_d$ with $\M_d=\{Q_{\t,d}:\t\in\SetR^d\}$
(say) can be reduced to a countable class $\M=\{P_d\}$ for which our
result holds, where $P_d$ is a Bayes or NML or other estimate of
$\M_d$ \cite{Gruenwald:07book}. Alternatively, $\bigcup_d\M_d$ could be
reduced to a countable class by considering only computable
parameters $\t$. Essentially all interesting model classes contain
such a countable topologically dense subset. Under certain
circumstances MDL still works for the non-computable parameters
\cite{Gruenwald:07book}. Alternatively one may simply reject
non-computable parameters on philosophical grounds
\cite{Hutter:04uaibook}. Finally, the techniques for the
countable case might aid proving general results for continuous
$\M$, possibly along the lines of \cite{Ryabko:09}.

\paradot{Contents}
The paper is organized as follows: In Section \ref{secFIP} we
provide some insights how MDL and Bayes work in restricted settings,
what breaks down for general countable $\M$, and how to circumvent
the problems. The formal development starts with Section
\ref{secNMR}, which introduces notation and our main result. The
proof for finite $\M$ is presented in Section \ref{secPrF} and for
denumerable $\M$ in Section \ref{secPrC}. In Section \ref{secAppl}
we show how the result can be applied to sequence prediction,
classification and regression, discriminative learning, and
reinforcement learning. Section \ref{secVar} discusses some MDL
variations.

\section{Facts, Insights, Problems}\label{secFIP}

Before starting with the formal development, we describe how MDL and
Bayes work in some restricted settings, what breaks down for general
countable $\M$, and how to circumvent the problems. For
deterministic environments, MDL reduces to learning by elimination,
and the four results in \req{results} can easily be understood.
Consistency of MDL for i.i.d.\ (and stationary-ergodic) sources is
also intelligible. For general $\M$, MDL may no longer converge to
the true model. We have to give up the idea of model identification,
and concentrate on predictive performance.

\paradot{Deterministic MDL = elimination learning}
For a countable class $\M=\{Q_1,Q_2,...\}$ of {\em deterministic}
theories=models=hypotheses=sequences, sorted w.r.t.\ to their
complexity=codelength $K(Q_i)=2\lb i$ (say) it is easy to see why
MDL works: Each $Q$ is a model for one infinite sequence
$x_{1:\infty}^Q$, i.e.\ $Q(x^Q)=1$. Given the true observations
$x\equiv x_{1:\l}^P$ so far, MDL selects the simplest $Q$ consistent
with $x_{1:\l}^P$ and for $h=1$ predicts $x_{\l+1}^Q$. This (and
potentially other) $Q$ becomes (forever) inconsistent if and only if
the prediction was wrong. Assume the true model is $\P=Q_m$. Since
elimination occurs in order of increasing index $i$, and $Q_m$ never
makes any error, MDL makes at most $m-1$ prediction errors. Indeed,
what we have described is just classical Gold style learning by
elimination. For $1<h<\infty$, the prediction $x_{\l+1:\l+h}^Q$ may
be wrong only on $x_{\l+h}^Q$, which causes $h$ wrong predictions
before the error is revealed. (Note that at time $\l$ only $x_\l^P$
is revealed.) Hence the total number of errors is bounded by $h\cdot
(m-1)$. The bound is for instance attained on the class consisting
of $Q_i= 1^{ih}0^\infty$, and the true sequence switches from 1 to 0
after having observed $m\cdot h$ ones. For $h=\infty$, a wrong
prediction gets {\em eventually} revealed. Hence each wrong $Q_i$
($i<m$) gets eventually eliminated, i.e.\ $P$ gets eventually
selected. So for $h=\infty$ we can (still/only) show that the number
of errors is finite. No bound on the {\em number} of errors in terms
of $m$ only is possible. For instance, for
$\M=\{Q_1=1^\infty,Q_2=P=1^n 0^\infty\}$, it takes $n$ time steps to
reveal that prediction $1^\infty$ is wrong, and $n$ can be chosen
arbitrarily large.

\paradot{Deterministic Bayes = majority learning}
Bayesian learning is at the same time, closely related to and very
different from MDL. Bayes predicts with a $w_Q$-weighted average of
the models (rather than with a single one). For a deterministic
class, Bayes is similar to prediction by majority: Consider the
models consistent with the true observation $x_{1:\l}^P$, having
total weight $W$, and take the weighted majority prediction (this is
the Bayes-optimal decision under 0-1 loss, Bayesian prediction would
randomize). For $h=1$, making a wrong prediction means that $Q$'s
contributing to at least half of the total weight $W$ get
eliminated. Since $P=Q_m$ never gets eliminated, we have $w_P\leq
W\leq 2^{-\#\text{Errors}}$, hence the number of errors is bounded
by $\lb w_P^{-1}$. For probabilistic Bayesian prediction proper, it
is also easy to see that the expected number of errors is bounded by
$\ln w_P^{-1}$. One can show that these bounds are essentially
sharp. (e.g. for $Q_i$ defined as the digits after the comma of the
binary expansion of $(i-1)/2^n$ for $i=1..m$ and $m=2^n-1$.) With
the same reasoning as in the MDL case, for $h>1$ we have to multiply
the bound by $h$; and for $h=\infty$ we get correct prediction
eventually, but no explicit bound anymore.

\paradot{Comparison of deterministic$\leftrightarrow$probabilistic and MDL$\leftrightarrow$Bayes}
The flavor of results carries over to some extent to the
probabilistic case. On a very abstract level even the line of
reasoning carries over, although this is deeply buried in the
sophisticated mathematical analysis of the latter. So the special
deterministic case illustrates the more complex probabilistic case.
For instance for $h=1$ and $w_i\propto 1/i^2$, we see that ``Bayes''
makes only $2\lb m$ errors, while MDL can make up to the $m$ errors.
This carries over to the probabilistic case. Also the multiplier $h$
for $1<h<\infty$ and the lack of an explicit bound for $h=\infty$
carries over. Cf.\ the bounds in \req{results}.
The reader is invited to reveal other relations not explicitly
mentioned here.
The differences are as follows: In the probabilistic case, the true
$P$ can in general not be identified anymore. Further, while the
Bayesian bound trivially follows from the 1/2-century old classical
merging of opinions result \cite{Blackwell:62}, the corresponding
MDL bound we prove in this paper is more difficult to
obtain.

\paradot{Consistency of MDL for stationary-ergodic sources}
For an i.i.d.\ class $\M$, the law of large numbers applied to the
random variables $Z_t:=\log[P(x_t)/Q(x_t)]$ implies
$\fr1\l\sum_{t=1}^\l Z_t\to\text{KL}(P||Q):=\sum_{x_1}P(x_1)\log
[P(x_1)/Q(x_1)]$ with $P$-probability 1. Either the Kullback-Leibler
(KL) divergence is zero, which is the case if and only if $P=Q$, or
$\log P(x_{1:\l})-\log Q(x_{1:\l}) \equiv \sum_{t=1}^\l Z_\l \sim
\text{KL}(P||Q)\l \to \infty$, i.e.\ asymptotically MDL does not
select $Q$. For countable $\M$, a refinement of this argument shows
that MDL eventually selects $P$ \cite{Barron:91}. This reasoning can
be extended to stationary-ergodic $\M$, but essentially not beyond.
To see where the limitation comes from, we present some troubling
examples.

\paradot{Trouble makers}
For instance, let $P$ be a Bernoulli$(\t_0)$ process, but let the
$Q$-probability that $x_t=1$ be $\t_t$, i.e.\ time-dependent (still
assuming independence). For a suitably converging but
``oscillating'' (i.e.\ infinitely often larger and smaller than its limit)
sequence $\t_t\to\t_0$ one can show that $\log[P(x_{1:t})/Q(x_{1:t})]$
converges to but oscillates around $K(Q)-K(P)$ w.p.1, i.e.\ there
are non-stationary distributions for which MDL does not converge
(not even to a wrong distribution).

One idea to solve this problem is to partition $\M$, where two
distributions are in the same partition if and only if they are
asymptotically indistinguishable (like $P$ and $Q$ above), and then
ask MDL to only identify a partition. This approach cannot succeed
generally, whatever particular criterion is used, for the following
reason:
Let $P(x_1)>0$ $\forall x_1$. For $x_1=1$, let $P$ and $Q$ be
asymptotically indistinguishable, e.g. $P=Q$ on the remainder of the
sequence. For $x_1=0$, let $P$ and $Q$ be asymptotically
distinguishable distributions, e.g.\ different Bernoullis. This
shows that for non-ergodic sources like this one, asymptotic
distinguishability depends on the drawn sequence. The first
observation can lead to totally different futures.

\paradot{Predictive MDL avoids trouble}
The Bayesian posterior does not need to converge to a single (true or
other) distribution, in order for prediction to work. We can do
something similar for MDL. At each time we still select a single
distribution, but give up the idea of identifying a single
distribution asymptotically. We just measure predictive success, and
accept infinite oscillations. That's the approach taken in this
paper.

\section{Notation and Main Result}\label{secNMR}

The formal development starts with this section. We need probability
measures and filters for infinite sequences, conditional
probabilities and densities, the total variation distance, and the
concept of merging (of opinions), in order to formally state our
main result.

\paradot{Measures on sequences}
Let $(\O:=\X^\infty,\F,\P)$ be the space
of infinite sequences with natural filtration and product
$\sigma$-field $\F$ and probability measure $\P$. Let $\o\in\O$ be
an infinite sequence sampled from the true measure $\P$. Except when
mentioned otherwise, all probability statements and expectations
refer to $\P$, e.g.\ almost surely (a.s.) and with probability 1 (w.p.1)
are short for with $P$-probability 1 (w.$P$.p.1).
Let $x=x_{1:\l}=\o_{1:\l}$ be the first $\l$ symbols of $\o$.

For countable $\X$, the probability that an infinite sequence starts
with $x$ is $P(x):=\P[\{x\}\times\X^\infty]$. The conditional
distribution of an event $A$ given $x$ is
$\P[A|x]:=\P[A\cap(\{x\}\times\X^\infty)]/P(x)$, which exists
w.p.1. For other probability measures $Q$ on $\O$, we define
$Q(x)$ and $Q[A|x]$ analogously. General $\X$ are considered at
the end of this section.

\paradot{Convergence in total variation}
$P$ is said to be {\em absolutely continuous} relative to $Q$, written
\beqn
  P\ll Q \quad:\Leftrightarrow\quad
  [Q[A]=0 \mbox{ implies } P[A]=0 \mbox{ for all } A\in\F]
\eeqn
$P$ and $Q$ are said to be {\em mutually singular}, written $P\bot Q$,
iff there exists an $A\in\F$ for which $P[A]=1$ and $Q[A]=0$.
The {\em total variation distance} (tvd) between
$Q$ and $P$ given $x$ is defined as
\beq\label{tvd}
  d(P,Q|x) \;:=\; \sup_{A\in\F}\big|Q[A|x]-P[A|x]\big|
\eeq
$Q$ is said to {\em predict $P$ in tvd} (or merge with $P$) if
$d(P,Q|x)\to 0$ for $\l(x)\to\infty$ with $P$-probability 1. Note
that this in particular implies, but is stronger than one-step
predictive on- and off-sequence convergence
$Q(x_{\l+1}=a_{\l+1}|x_{1:\l})-P(x_{\l+1}=a_{\l+1}|x_{1:\l})\to 0$ for any
$a$, not necessarily equal $\o$ \cite{Kalai:94}.
The famous Blackwell and Dubins convergence result
\cite{Blackwell:62} states that if $P$ is absolutely continuous
relative to $Q$, then (and only then \cite{Kalai:94}) $Q$ merges with
$P$:
\beqn
  \qmbox{If} P\ll Q \qmbox{then} d(P,Q|x)\to 0 \qmbox{w.p.1 \ for} \l(x)\to\infty
\eeqn

\paradot{Bayesian prediction}
This result can immediately be utilized for Bayesian prediction. Let
$\M:=\{Q_1,Q_2,Q_3,...\}$ be a countable (finite or infinite) class
of probability measures, and $\Bayes[A]:=\sum_{Q\in\M} Q[A]w_Q$ with
$w_Q>0$ $\forall Q$ and $\sum_{Q\in\M}w_Q=1$. If the model
assumption $P\in\M$ holds, then obviously $P\ll\Bayes$, hence $\Bayes$ merges
with $P$, i.e.\ $d(P,\Bayes|x)\to 0$ w.p.1 for all $P\in\M$. Unlike many
other Bayesian convergence and consistency theorems, no
(independence, ergodicity, stationarity, identifiability, or
other) assumption on the model class $\M$ need to be made. Good
convergence rates for the weaker $d_{h<\infty}$ distances have also
been shown \cite{Hutter:03spupper}.
The analogous result for MDL is as follows:

\begin{theorem}[MDL predictions]\label{thmMDLtvp}
Let $\M$ be a countable class of probability measures on $\X^\infty$
containing the unknown true sampling distribution $P$. No
(independence, ergodicity, stationarity, identifiability, or
other) assumptions need to be made on $\M$. Let
\beqn
  \MDL^x \;:=\; \arg\min_{Q\in\M}\{-\log Q(x) +K(Q)\}
  \qmbox{with} \sum_{Q\in\M}2^{-K(Q)}<\infty
\eeqn
be the measure selected by MDL at time $\l$ given $x\in\X^\l$.
Then the predictive distributions $\MDL^x[\cdot|x]$ converge to
$P[\cdot|x]$ in the sense that
\beqn
  d(P,\MDL^x|x) \;\equiv\; \sup_{A\in\F}\big|\MDL^x[A|x]-P[A|x]\big|
  \;\to\; 0 \qmbox{for} \l(x)\to\infty \qmbox{w.p.1}
\eeqn
\end{theorem}
$K(Q)$ is usually interpreted and defined as the length of some
prefix code for $Q$, in which case $\sum_Q 2^{-K(Q)}\leq 1$. If
$K(Q):=\lb w_Q^{-1}$ is chosen as complexity, by Bayes rule
$\Prob(Q|x)=Q(x)w_Q/\Bayes(x)$, the maximum a posteriori estimate
$\text{MAP}^x := \arg\max_{Q\in\M}\{\Pr(Q|x)\} \equiv \MDL^x$. Hence
the theorem also applies to MAP. The proof of the theorem is
surprisingly subtle and complex compared to the analogous Bayesian
case. One reason is that $\MDL^x(x)$ is not a measure on
$\X^\infty$.

\paradot{Arbitrary $\X$}
For arbitrary $\X$, definitions are more subtle.
The casual reader satisfied with countable $\X$ can skip this paragraph.
We can consider even more generally $x_t\in\X_t$
\cite{Blackwell:62}.
Let $\B_t$ be a $\sigma$-field of subsets of $\X_t$ for
$t=1,2,3,...$. Let $\F_\l$ be the $\sigma$-field for
$\X^\l:=\X_1\times...\times\X_\l$ generated by (i.e.\ the smallest
$\sigma$-field containing) $\B_1\times...\times\B_\l$ for
$\l\leq\infty$. Let $(\O:=\X^\infty,\F=\F_\infty,P)$ be a
probability space. Let $P_\l$ be the marginal distribution on
$(\X^\l,\F_\l)$, i.e.\
$P_\l[A]:=P[A\times\X_{\l+1}\times\X_{\l+2}\times...]$ for
$A\in\F_\l$. The predictive distribution $P^\l[A|x_{1:\l}]$ is (a
version of) the conditional distribution of the future
``$x_{\l+1:\infty}$'' given past $x_{1:\l}$, implicitly defined by
$\int P^\l[A|x_{1:\l}] dP_\l(x_{1:\l}):=P[A]$ $\forall A\in\F$.
Similarly define $Q_\l$ and $Q^\l$ for the other $Q\in\M$. See
\cite{Doob:53} for details.

Let $M$ be a measure on $\O$ such that $Q$ is absolutely continuous
(see below) relative to\ $M$ for all $Q\in\M$. For instance
$M[\cdot]=\Bayes[\cdot]$ has this property. Now define the density
(Radon-Nikodym derivative) $Q_\l(x_{1:\l})$ (round brackets) of
measure $Q_\l[\cdot]$ (square brackets) relative to $M_\l[\cdot]$.
It is important to note that all essential quantities, in particular
$\MDL^x$, are independent of the particular choice of $M$. We
therefore plainly speak of the $Q$-density or even $Q$-probability
of $x$.

For countable $\X$ and counting measure $M$, $Q^\l[A|x]$ and
$Q_\l(x)$ coincide with $Q[A|x]$ and $Q(x)$ above. In the
following, we drop the sup\&\-superscripts $\l$, since they will
always be clear from the argument.
Note that by Carathéodory's extension theorem, $\{Q(x):x\in\X^*\}$
uniquely defines $Q[A]$ $\forall A\in\F$.

\section{Proof for Finite Model Class}\label{secPrF}

We first prove Theorem \ref{thmMDLtvp} for finite model classes
$\M$. For this we need the following Definition and Lemma:

\begin{definition}[Relations between $Q$ and $P$]\label{defRBM}
For any probability measures $Q$ and $P$, let
\begin{itemize}\parskip=0ex\parsep=0ex\itemsep=0ex
\item $Q^r+Q^s=Q$ be the
Lebesgue decomposition of $Q$ relative to $P$ into an absolutely
continuous non-negative measure $Q^r\ll P$ and a singular
non-negative measure $Q^s\bot P$.
\item $g(\o) \;:=\; dQ^r/dP \;=\; \lim_{\l\to\infty}[Q(x_{1:\l})/P(x_{1:\l})]$ be (a
version of) the Radon-Nikodym derivative, i.e.\ $Q^r[A]=\int_A
g\,dP$.
\item $\Oo   :=\; \{\o: Q(x_{1:\l})/P(x_{1:\l})\to 0\} \;\equiv\; \{\o:g(\o)=0\}$.
\item $\Ov \;:=\; \{\o : d(P,Q|x)\to 0 \mbox{ for } \l(x)\to\infty \}$.
\end{itemize}
\end{definition}
It is well-known that the Lebesgue decomposition exists and is
unique. The representation of the Radon-Nikodym derivative as a
limit of local densities can e.g.\ be found in \cite[VII\S
8]{Doob:53}: $Z_\l^{r/s}(\o):=Q^{r/s}(x_{1:\l})/P(x_{1:\l})$ for
$\l=1,2,3,...$ constitute two martingale sequences, which converge
w.p.1. $Q^r\ll P$ implies that the limit $Z_\infty^r$ is the
Radon-Nikodym derivative $dQ^r/dP$. (Indeed, Doob's martingale
convergence theorem can be used to prove the Radon-Nikodym theorem.)
$Q^s\bot P$ implies $Z_\infty^r=0$ w.p.1. So $g$ is uniquely defined
and finite w.p.1.

\begin{lemma}[Generalized merging of opinions]\label{lemGMO}
For any $Q$ and $P$, the following holds:
\begin{itemize}\parskip=0ex\parsep=0ex\itemsep=0ex
\item[(i)] $P\ll Q$ if and only if $\P[\Oo]=0$
\item[(ii)] $\P[\Oo]=0$ implies $\P[\Ov]=1$ \hfill [(i)+\cite{Blackwell:62}]
\item[(iii)] $\P[\Oo\cup\Ov]=1$ \hfill [generalizes (ii)]
\end{itemize}
\end{lemma}
$(i)$ says that $Q(x)/P(x)$ converges almost surely to a strictly
positive value if and only if $P$ is absolutely continuous relative to
$Q$,
$(ii)$ says that an almost sure positive limit of $Q(x)/P(x)$ implies
that $Q$ merges with $P$.
$(iii)$ says that even if $P\not\ll Q$, we still have $d(P,Q|x)\to 0$
on almost every sequence that has a positive limit of $Q(x)/P(x)$.

\paradot{Proof}
Recall Definition \ref{defRBM}.

$(i\!\Leftarrow)$ Assume $P[\Oo]=0$: $P[A]>0$ implies $Q[A]\geq
Q^r[A]=\int_A g\,dP>0$, since $g>0$ a.s.\ by assumption $P[\Oo]=0$.
Therefore $P\ll Q$.

$(i\!\!\Rightarrow)$ Assume $P\ll Q$: Choose a $B$ for which $P[B]=1$
and $Q^s[B]=0$. Now $Q^r[\Oo]=\int_\Oo g\,dP=0$ implies $0\leq
Q[B\cap\Oo]\leq Q^s[B]+Q^r[\Oo] = 0+0$. By $P\ll Q$ this implies
$P[B\cap\Oo]=0$, hence $P[\Oo]=0$.

$(ii)$ That $P\ll Q$ implies $\P[\Ov]=1$ is Blackwell-Dubins'
celebrated result. The result now follows from (i).

$(iii)$ generalizes \cite{Blackwell:62}. For $\P[\Oo]=0$ it reduces
to (ii). The case $\P[\Oo]=1$ is trivial. Therefore we can assume
$0<\P[\Oo]<1$. Consider measure $P'[A]:=\P[A|B]$ conditioned on
$B:=\O\setminus\Oo$.

Assume $Q[A]=0$. Using $\int_{\Oo}g\,dP=0$, we get $0=Q^r[A]=\int_A
g\,dP=\int_{A\setminus\Oo}g\,dP$. Since $g>0$ outside $\Oo$, this
implies $P[A\setminus\Oo]=0$. So $P'[A]=P[A\cap
B]/P[B]=P[A\setminus\Oo]/P[B]=0$. Hence $P'\ll Q$. Now (ii) implies
$d(P',Q|x)\to 0$ with $P'$ probability 1. Since $P'\ll P$ we also
get $d(P',P|x)\to 0$ w.$P'$.p.1.

Together this implies $0\leq d(P,Q|x)\leq d(P',P|x)+d(P',Q|x)\to 0$
w.$P'$.p.1, i.e.\ $P'[\Ov]=1$. The claim now follows from
\hfill\raisebox{-8ex}[0ex][0ex]{\qed}
\bqan
  P[\Oo\cup\Ov]
  &=& P'[\Oo\cup\Ov]P[\O\setminus\Oo] + P[\Oo\cup\Ov|\Oo]P[\Oo] \\
  &=& 1\cdot P[\O\setminus\Oo] + 1\cdot P[\Oo] \;=\; P[\O] \;=\; 1
\eqan

The intuition behind the proof of Theorem \ref{thmMDLtvp} is as
follows. MDL will asymptotically not select $Q$ for which
$Q(x)/P(x)\to 0$. Hence for those $Q$ potentially selected by MDL,
we have $\o\not\in\Oo$, hence $\o\in\Ov$, for which $d(P,Q|x)\to 0$
(a.s.). The technical difficulties are for finite $\M$ that the
eligible $Q$ depend on the sequence $\o$, and for infinite $\M$ to
deal with non-uniformly converging $d$, i.e.\ to infer
$d(P,\MDL^x|x)\to 0$.

\paradot{Proof of Theorem \ref{thmMDLtvp} for finite $\M$}
Recall Definition \ref{defRBM}, and let
$g_Q,\Oo_{\!\!\!\!\!\;Q},\Ov_Q$ refer to some
$Q\in\M\equiv\{Q_1,...,Q_m\}$. The set of sequences $\o$ for which
some $g_Q$ for some $Q\in\M$ is undefined has $P$-measure zero, and
hence can be ignored. Fix some sequence $\o\in\O$ for which
$g_Q(\o)$ is defined for all $Q\in\M$, and let
$\M_\o:=\{Q\in\M:g_Q(\o)=0\}$.
\beqn
  \MDL^x:=\arg\min_{Q\in\M}L_Q(x), \qmbox{where}
  L_Q(x):=-\log Q(x)+K(Q).
\eeqn
Consider the difference
\beqn
  L_Q(x)-L_P(x) \;=\; -\log{Q(x)\over P(x)}+K(Q)-K(P)
  \;\toinfty{\l} -\log g_Q(\o) +K(Q)-K(P)
\eeqn
For $Q\in\M_\o$, the r.h.s.\ is $+\infty$, hence
\beqn
  \forall Q\!\in\!\M_\o\,\exists \l_Q \forall \l\!>\!\l_Q : L_Q(x)>L_P(x)
\eeqn
Since $\M$ is finite, this implies
\beqn
  \forall \l\!>\!\l_0\,\forall Q\!\in\!\M_\o : L_Q(x)>L_P(x),
  \qmbox{where} \l_0:=\max\{\l_Q:Q\in\M_\o\} < \infty
\eeqn
Therefore, since $P\in\M$, we have $\MDL^x\not\in\M_\o$ $\forall
\l>\l_0$, so we can safely ignore all $Q\in\M_\o$ and
focus on $Q\in\overline\M_\o:=\M\setminus\M_\o$. Let
$\O_1:=\bigcap_{Q\in\overline\M_\o}(\Oo_{\!\!\!\!\!\;Q}\cup\Ov_Q)$.
Since $P[\O_1]=1$ by Lemma \ref{lemGMO}(iii), we can also assume
$\o\in\O_1$.
\beqn
  Q\in\overline\M_\o \quad\Rightarrow\quad g_Q(\o)>0
  \quad\Rightarrow\quad \o\not\in\Oo_{\!\!\!\!\!\;Q}
  \quad\Rightarrow\quad \o\in\Ov_Q
  \quad\Rightarrow\quad d(P,Q|x)\to 0
\eeqn
This implies
\vspace{-2ex}\beqn
  d(\P,\MDL^x|x) \;\leq\; \sup_{Q\in\overline\M_\o} d(\P,Q|x) \;\to\; 0
\eeqn
where the inequality holds for $\l>\l_0$ and the limit holds, since
$\M$ is finite. Since the set of $\o$ excluded in our considerations
has measure zero, $d(\P,\MDL^x|x)\to 0$ w.p.1, which proves the
theorem for finite $\M$.
\qed

\section{Proof for Countable Model Class}\label{secPrC}

The proof in the previous Section crucially exploited finiteness of
$\M$. We want to prove that the probability that MDL asymptotically
selects ``complex'' $Q$ is small. The following Lemma establishes
that the probability that MDL selects a {\em specific} complex $Q$
infinitely often is small.

\begin{lemma}[MDL avoids complex probability measures $Q$]\label{lemNCQ}
For any $Q$ and $P$ we have
$\P[Q(x)/P(x)\geq c \mbox{ infinitly often}]\leq 1/c$.
\end{lemma}

\paradot{Proof}
\vspace{-4ex}\bqan
  & & \P[\forall \l_0\exists \l\!>\!\l_0\!:\!\frac{Q(x)}{P(x)}\geq c]
  \;\stackrel{(a)}=\; \P[\mathop{\overline\lim}_{\l\to\infty} \frac{Q(x)}{P(x)}\geq c] \;\leq\;
\\
  & & \;\stackrel{(b)}\leq\; \frac1c\E[\mathop{\overline\lim}_\l \frac{Q(x)}{P(x)}]
  \;\stackrel{(c)}=\; \frac1c\E[\mathop{\underline\lim}_\l \frac{Q(x)}{P(x)}]
  \;\stackrel{(d)}\leq\; \frac1c\mathop{\underline\lim}_\l\E[\frac{Q(x)}{P(x)}]
  \;\stackrel{(e)}\equiv\; \frac1c
\eqan
(a) is true by definition of the limit superior $\overline\lim$, (b)
is Markov's inequality, (c) exploits the fact that the limit of
$Q(x)/P(x)$ exists w.p.1, (d) uses Fatou's lemma, and (e) is
obvious. \qed

For sufficiently complex $Q$, Lemma \ref{lemNCQ} implies that
$L_Q(x) > L_P(x)$ for most $x$. Since convergence is non-uniform in
$Q$, we cannot apply the Lemma to all (infinitely many) complex $Q$
directly, but need to lump them into one $\bar Q$.

\paradot{Proof of Theorem \ref{thmMDLtvp} for countable $\M$}
Let the $Q\in\M=\{Q_1,Q_2,...\}$ be ordered somehow, e.g.\ in
increasing order of complexity $K(Q)$, and $P=Q_n$. Choose some
(large) $m\geq n$ and let $\widetilde\M:=\{Q_{m+1},Q_{m+2},...\}$ be
the set of ``complex'' $Q$. We show that the probability that MDL
selects infinitely often complex $Q$ is small:
\bqan
  & & \nq\nq \P[\MDL^x\in\widetilde\M\mbox{ infinitely often}] \\
  &\equiv& \P[\forall \l_0\exists \l\!>\!\l_0:\MDL^x\in\widetilde\M] \\
  &\leq& \P[\forall \l_0\exists \l\!>\!\l_0\wedge Q\in\widetilde\M: L_Q(x)\leq L_P(x)] \\
  &=& \P[\forall \l_0\exists \l\!>\!\l_0 : \mathop{\smash\sup}_{i>m}\textstyle{Q_i(x)\over P(x)}2^{K(P)-K(Q_i)}\geq 1] \\
  &\smash{\stackrel{(a)}\leq}& \P[\forall \l_0\exists \l\!>\!\l_0 : \textstyle{\bar Q(x)\over P(x)}\,\delta\, 2^{K(P)}\geq 1] \\
  &\stackrel{(b)}\leq& \delta\, 2^{K(P)}
  \;\stackrel{(c)}\leq\; \eps
\eqan
The first three relations follow immediately from the definition of
the various quantities. Bound (a) is the crucial ``lumping'' step.
First we bound
\beqn
  \sup_{i>m}{Q_i(x)\over P(x)}2^{-K(Q_i)}
  \;\leq\; \sum_{i=m+1}^\infty {Q_i(x)\over P(x)}2^{-K(Q_i)}
  \;=\; \delta\, {\bar Q(x)\over P(x)},
\eeqn
\beqn
  \delta:=\sum_{i>m}2^{-K(Q_i)}<\infty,\qquad
  \bar Q(x) :={1\over\delta}\sum_{i>m}Q_i(x)2^{-K(Q_i)},
\eeqn
While $\MDL^{\displaystyle\cdot}[\cdot]$ is not a (single) measure
on $\O$ and hence difficult to deal with, $\bar Q$ is a proper
probability measure on $\Omega$. In a sense, this step reduces MDL
to Bayes. Now we apply Lemma \ref{lemNCQ} in (b) to the (single)
measure $\bar Q$. The bound (c) holds for sufficiently large $m=m_\eps(P)$,
since $\delta\to 0$ for $m\to\infty$.
This shows that for the sequence of MDL estimates
\beqn
  \{\MDL^{x_{1:\l}}\!:\!\l>\l_0\} \subseteq \{Q_1,...,Q_m\}
  \qmbox{with probability at least} 1-\eps
\eeqn
Hence the already proven Theorem \ref{thmMDLtvp} for finite
$\M$ implies that $d(\P,\MDL^x|x)\to 0$ with probability at least
$1-\eps$. Since convergence holds for every $\eps>0$, it holds w.p.1.
\qed

\section{Implications}\label{secAppl}

Due to its generality, Theorem \ref{thmMDLtvp} can be applied to
many problem classes. We illustrate some immediate implications of
Theorem \ref{thmMDLtvp} for time-series forecasting, classification,
regression, discriminative learning, and reinforcement learning.

\paradot{Time-series forecasting}
Classical online sequence prediction is concerned with predicting
$x_{\l+1}$ from (non-i.i.d.) sequence $x_{1:\l}$ for $\l=1,2,3,...$.
Forecasting farther into the future is possible by predicting
$x_{\l+1:\l+h}$ for some $h>0$. One can show that
$0\leq d_1 \leq d_h \leq d_{h+1} \leq d_\infty = 2d\leq 2$, see
\req{dh} and \req{tvd}. Hence Theorem \ref{thmMDLtvp} implies
good asymptotic (multi-step) predictions.
Offline learning is concerned with training a predictor on
$x_{1:\l}$ for fixed $\l$ in-house, and then selling and using the
predictor on $x_{\l+1:\infty}$ without further learning. Theorem
\ref{thmMDLtvp} shows that for enough training data, predictions
``post-learning'' will be good.

\paradot{Classification and Regression}
In classification (discrete $\X$) and regression (continuous $\X$),
a sample is a set of pairs $D=\{(y_1,x_1),...,(y_\l,x_\l)\}$, and a
functional relationship $\dot x=f(\dot y)$+noise, i.e.\ a
conditional probability $P(\dot x|\dot y)$ shall be learned. For
reasons apparent below, we have swapped the usual role of $\dot x$
and $\dot y$. The dots indicate $\dot x\in\X$ and $\dot y\in\Y)$, while
$x=x_{1:\l}\in\X^\l$ and $y=y_{1:\l}\in\Y^\l$.
%
If we assume that also $\dot y$ follows some distribution, and start
with a countable model class $\M$ of joint distributions $Q(\dot
x,\dot y)$ which contains the true joint distribution $P(\dot x,\dot
y)$, our main result implies that $\MDL^D[(\dot x,\dot y)|D]$
converges to the true distribution $P(\dot x,\dot y)$. Indeed
since/if samples are assumed i.i.d., we don't need to invoke our
general result.

\paradot{Discriminative learning}
Instead of learning a generative \cite{Jebara:03} joint
distribution $P(\dot x,\dot y)$, which requires model assumptions on
the input $\dot y$, we can discriminatively \cite{Long:07} learn
$P(\cdot|\dot y)$ directly without any assumption on $y$ (not even
i.i.d). We can simply treat $y_{1:\infty}$ as an oracle to all $Q$,
define $\M'=\{Q'\}$ with $Q'(x):=Q(x|y_{1:\infty})$, and apply our
main result to $\M'$, leading to $\MDL'{}^x[A|x]\to P'[A|x]$, i.e.\
$\MDL^{x|y_{1:\infty}}[A|x,y_{1:\infty}]\to P[A|x,y_{1:\infty}]$.
This not yet useful since $y_{1:\infty}$ is never known completely.
If $x_1,x_2,...$ are conditionally independent, we can write
\beqn
  Q(x|y) \;=\; \prod_{t=1}^\l Q(x_t|y_t)
  \;=\; \sum_{x_{\l+1:m}} \prod_{t=1}^m Q(x_t|y_t)
  \;=\; \sum_{x_{\l+1:m}} Q(x_{1:m}|y_{1:m})
  \;=\; Q(x_{1:\l}|y_{1:m})
\eeqn
Taking the limit $m\to\infty$ we get $Q(x|y)=Q(x|y_{1:\infty})$.
This is a generic property satisfied for all {\em causal} processes,
that a future $y_t$ for $t>l$ does not influence past observations
$x_{1:\l}$. Hence for a class of conditionally independent
distributions, we get $\MDL^{x|y}[A|x,y]\to P[A|x,y]$.
Since the $x$ given $y$ are {\em not} identically distributed,
classical MDL consistency results for i.i.d.\ or stationary-ergodic
sources do {\em not} apply. The following corollary formalizes our
findings:

\begin{corollary}[Discriminative MDL]\label{lemDMDL}
Let $\M\ni P$ be a class of discriminative causal distributions
$Q[\cdot|y_{1:\infty}]$, i.e.\ $Q(x|y_{1:\infty})=Q(x|y)$, where
$x=x_{1:\l}$ and $y=y_{1:\l}$. Regression and classification are
typical examples. Further assume $\M$ is countable. Let
$\MDL^{x|y}:=\arg\min_{Q\in\M}\{-\log Q(x|y)+K(Q)\}$ be the
discriminative MDL measure (at time $\l$ given $x,y$). Then
$\sup_A\big|\MDL^{x|y}[A|x,y]-P[A|x,y]\big|\to 0$ for
$\l(x)\to\infty$, $P[\cdot|y_{1:\infty}]$ almost surely, for {\em
every} sequence $y_{1:\infty}$.
\end{corollary}

For finite $\Y$ and conditionally independent $x$, the intuitive
reason how this can work is as follows: If $\dot y$ appears in
$y_{1:\infty}$ only finitely often, it plays asymptotically no role;
if it appears infinitely often, then $P(\cdot|\dot y)$ can be
learned. For infinite $\Y$ and deterministic $\M$, the result is
also intelligible: Every $\dot y$ might appear only once, but
probing enough function values $x_t=f(y_t)$ allows to identify the
function.

\paradot{Reinforcement learning (RL)}
In the agent framework \cite{Russell:03}, an agent interacts with an
environment in cycles. At time $t$, an agent chooses an action $y_t$
based on past experience $x_{<t}\equiv(x_1,...,x_{t-1})$ and past
actions $y_{<t}$ with probability $\pi(y_t|x_{<t}y_{<t})$ (say). This
leads to a new perception $x_t$ with probability
$\mu(x_t|x_{<t}y_{1:t})$ (say). Then cycle $t+1$ starts. Let
$P(xy)=\prod_{t=1}^\l \mu(x_t|x_{<t}y_{1:t})\pi(y_t|x_{<t}y_{<t})$
be the joint interaction probability. We make no (Markov,
stationarity, ergodicity) assumption on $\mu$ and $\pi$. They may be
POMDPs or beyond.

\begin{corollary}[Single-agent MDL]\label{lemSAMDL}
For a fixed policy=agent $\pi$, and a class of environments
$\{\nu_1,\nu_2,...\}\ni\mu$, let $\M=\{Q_i\}$ with
$Q_i(x|y)=\prod_{t=1}^\l \nu_i(x_t|x_{<t}y_{1:t})$. Then
$d(P[\cdot|y],\MDL^{x|y}) \to 0$ with joint $P$-probability 1.
\end{corollary}

The corollary follows immediately from the previous corollary and
the facts that the $Q_i$ are causal and that with
$P[\cdot|y_{1:\infty}]$-probability 1 $\forall y_{1:\infty}$ implies
w.$P$.p.1 jointly in $x$ and $y$.

In reinforcement learning \cite{Sutton:98}, the perception
$x_t:=(o_t,r_t)$ consists of some regular observation $o_t$ and
a reward $r_t\in[0,1]$. Goal is to find a policy which maximizes
accrued reward in the long run. The previous corollary implies

\begin{corollary}[Fixed-policy MDL value function convergence]\label{lemFPMDL}
Let $V_P[xy]:=\E_{P[\cdot|xy]}[r_{\l+1}+\g r_{\l+2} + \g^2 r_{\l+3}
+ ...]$ be the future $\g$-discounted $P$-expected reward sum (true value
of $\pi$), and similarly $V_{Q_i}[xy]$ for $Q_i$. Then the \MDL\
value converges to the true value, i.e.\
$V_{\MDL^{x|y}}[xy]-V_P[xy]\to 0$, w.$P$.p.1.\ for any policy $\pi$.
\end{corollary}

\paradot{Proof}
The corollary follows from the general inequality
\beqn
  \big|\E_P[f]-\E_Q[f]\big| \;\leq\; \sup|f|\cdot \sup_A\big|P[A]-Q[A]\big|
\eeqn
by inserting $f:=r_{\l+1}+\g r_{\l+2} + \g^2 r_{\l+3} + ...$ and
$P=P[\cdot|xy]$ and $Q=\MDL^{x|y}[\cdot|xy]$, and using $0\leq f\leq
1/(1-\g)$ and Corollary \ref{lemSAMDL}.
\qed

Since the value function probes the infinite future, we really made
use of our convergence result in total variation. Corollary
\ref{lemFPMDL} shows that MDL approximates the true value
asymptotically arbitrarily well. The result is weaker than it may
appear. Following the policy that maximizes the estimated (MDL)
value is often not a good idea, since the policy does not explore
properly \cite{Hutter:04uaibook}. Nevertheless, it is a
reassuring non-trivial result.

\section{Variations}\label{secVar}

MDL is more a general principle for model selection than a uniquely
defined procedure. For instance, there are {\em crude} and {\em
refined} MDL \cite{Gruenwald:07book}, the related {\em MML}
principle \cite{Wallace:05}, a {\em static}, a {\em dynamic}, and a
{\em hybrid} way of using MDL for prediction \cite{Hutter:05mdl2px},
and other variations.
For our setup, we could have defined multi-step lookahead prediction as a
product of single-step predictions:
\beqn
  \MDLI(x_{1:\l}) \;:=\; \prod_{t=1}^\l \MDL^{x_{<t}}(x_t|x_{<t}), \qquad
  \MDLI(z|x) \;=\; \MDLI(xz)/\MDLI(x)
\eeqn
which is a more {\em incremental} MDL version. Both, $\MDL^x$ and
$\MDLI$ are `static' in the sense of \cite{Hutter:05mdl2px}, and
each allows for a dynamic and a hybrid version.
Due to its incremental nature, $\MDLI$ likely has better predictive
properties than $\MDL^x$, and conveniently defines a {\em single}
measure over $\X^\infty$, but inconveniently is $\not\in\M$.
One reason for using $\MDL$ is that it can be computationally
simpler than Bayes. E.g. if $\M$ is a class of MDPs, then $\MDL^x$
is still an MDP and hence tractable, but $\MDLI$ like Bayes are a
nightmare to deal with.

\paradot{Acknowledgements}
My thanks go to Peter Sunehag for useful discussions.

{\addcontentsline{toc}{section}{\refname}
\begin{small}

\end{small}

\end{document}